\documentclass[11pt,leqno]{amsart}

\usepackage{amssymb}
\usepackage{amsmath}
\usepackage{amsthm}
\usepackage{amsfonts}
\usepackage{bm,bbm}
\usepackage{mathrsfs}
\usepackage{centernot}
\usepackage{pgfplots}
\usepackage{fancybox}
\usepackage{enumitem}

\usepackage{a4wide}
\usepackage{bookmark}
\usepackage{hyperref}
\hypersetup{pdfstartview={FitH}}

\newcommand{\n}{_\textup{left}}
\newcommand{\m}{_\textup{right}}

\newtheorem{theorem}{Theorem}

\newcommand{\R}{\mathbb{R}}

\newcommand{\B}{\mathcal{B}}
\newcommand{\A}{\mathcal{A}}
\newcommand{\F}{\mathcal{F}}
\newcommand{\E}{\mathbb{E}}

\newcommand{\D}{\mathcal{D}}

\numberwithin{equation}{section}

\begin{document}

\title[Bellman functions for paraproducts]{Bellman functions and $\textup{L}^p$ estimates for paraproducts}

\author[V. Kova\v{c}]{Vjekoslav Kova\v{c}}
\address{Vjekoslav Kova\v{c}, Department of Mathematics, Faculty of Science, University of Zagreb, Bijeni\v{c}ka cesta 30, 10000 Zagreb, Croatia}
\email{vjekovac@math.hr}

\author[K. A. \v{S}kreb]{Kristina Ana \v{S}kreb}
\address{Kristina Ana \v{S}kreb, Faculty of Civil Engineering, University of Zagreb, Fra Andrije Ka\v{c}i\'{c}a Mio\v{s}i\'{c}a 26, 10000 Zagreb, Croatia}
\email{kskreb@grad.hr}


\subjclass[2010]{Primary 60G46; 
Secondary 42B15} 


\begin{abstract}
We give an explicit formula for one possible Bellman function associated with the $\textup{L}^p$ boundedness of dyadic paraproducts regarded as bilinear operators or trilinear forms. Then we apply the same Bellman function in various other settings, to give self-contained alternative proofs of the estimates for several classical operators. These include the martingale paraproducts of Ba\~{n}uelos and Bennett and the paraproducts with respect to the heat flows.
\end{abstract}

\maketitle


\section{Introduction}
According to Janson and Peetre \cite{JP:par} the name ``paraproduct'' denotes an idea rather than a unique object. Various types of paraproducts appear in the literature on analysis or probability and in each case certain boundedness properties (i.e.\@ continuity) are crucial for their applications. An interested reader can find the historical overview and further references in the short expository paper \cite{BMN}. In this paper we will focus mostly on martingale paraproducts and revisit the $\textup{L}^p$ estimates, which they are well-known to satisfy.

We start with the dyadic paraproduct as a motivation for the forthcoming Bellman function that we construct. For two functions $f$ and $g$ from an appropriate space of real-valued test functions on $\mathbb{R}$ we can define the \emph{dyadic paraproduct} as a bilinear operator in the following way:
\begin{equation}
\Pi_{\epsilon}(f,g) := \sum_{I\in\mathcal{D}} \epsilon_I  |I|^{-2} \langle f,\mathbbm{1}_I\rangle \langle g,\mathbbm{h}_I \rangle \mathbbm{h}_I.
\label{dydpar:bil}
\end{equation}
Here $\D$ denotes the family of dyadic intervals in $\R$, $\mathbbm{1}_I$ is the indicator function of an interval $I$, $\mathbbm{h}_I:=\mathbbm{1}_{I\n}-\mathbbm{1}_{I\m}$ is the $\textup{L}^\infty$-normalized Haar function, while $I\n$ and $I\m$ are respectively the left half and the right half of $I$. Moreover, $\langle\cdot,\cdot\rangle$ denotes the standard inner product with respect to the Lebesgue measure and $\epsilon=(\epsilon_I)_{I\in\D}$ is a collection of real numbers such that $|\epsilon_I|\leqslant 1$ for each $I\in\mathcal{D}$. (If we choose $\epsilon_I \in \{-1,1\}$, then they simply represent $-$ and $+$ signs.) A convenient choice for the test functions are the so-called \emph{dyadic step functions}, i.e.\@ finite linear combinations of the indicator functions of dyadic intervals.

Typically, such an object is viewed as a linear operator in $g$ with $f$ fixed, when it becomes a particular instance of \emph{Bukholder's martingale transform} \cite{Burk1}. Alternatively, one can fix $g$ and consider it as a linear operator in $f$, in which case it is known as the linear paraproduct. In this text we prefer to look at $\Pi_{\epsilon}$ symmetrically and discuss its properties as a bilinear operator. This is partly motivated by the multilinear harmonic analysis, where more singular operators of this type are studied; see the book \cite{Thie}.

Equivalently, we can define the dyadic paraproduct as a trilinear form. We take the third test function $h$, and dualize \eqref{dydpar:bil} to get
\begin{align} \label{dydpar:tri}
\Lambda_{\epsilon}(f,g,h):=& \int_\R \Pi_{\epsilon}(f,g) h = \sum_{I\in \D}\epsilon_I |I|^{-2}\langle f, \mathbbm{1}_I\rangle \langle g, \mathbbm{h}_I\rangle \langle h,\mathbbm{h}_I \rangle \\
=&\sum_{I\in\mathcal{D}} \epsilon_I |I|[f]_I\Big( \frac{[g]_{I\n}-[g]_{I\m}}{2}\Big)\Big( \frac{[h]_{I\n}-[h]_{I\m}}{2}\Big). \nonumber
\end{align}
Here $[f]_I$ denotes the average of a function $f$ on a dyadic interval $I$.

It is well known that \eqref{dydpar:tri} satisfies certain $\textup{L}^p$ estimates, i.e.\@ there exists a finite constant $\mathcal{C}_{p,q,r}\geqslant 0$ depending only on three exponents $p,q,r$ such that
\begin{equation}
|\Lambda_{\epsilon}(f,g,h)| \leqslant \mathcal{C}_{p,q,r} \|f\|_{\textup{L}^p(\R)}\|g\|_{\textup{L}^q(\R)}\|h\|_{\textup{L}^r(\R)}
\label{lp:trilinear}
\end{equation}
holds whenever $1<p,q,r\leqslant \infty$ and $\frac{1}{p}+\frac{1}{q}+\frac{1}{r}=1$. By $\|\cdot\|_{\textup{L}^p(\R)}$ we have denoted the $\textup{L}^p$ norm on $\R$ with respect to the Lebesgue measure.

The easiest proof of \eqref{lp:trilinear} when $q,r<\infty$ uses boundedness of the dyadic maximal function and the dyadic square function. We simply apply the Cauchy-Schwarz and H\"{o}lder's inequality to get
\[|\Lambda_{\epsilon}(f,g,h)| \leqslant \int_\R (Mf)(Sg)(Sh) \leqslant \|Mf\|_{\textup{L}^p(\R)}\|Sg\|_{\textup{L}^q(\R)}\|Sh\|_{\textup{L}^r(\R)},\]
where
\[ Mf:=\sup_{I\in \D}|I|^{-1} |\langle f, \mathbbm{1}_I\rangle| \mathbbm{1}_I \quad \text{and} \quad Sf:=\Big(\sum_{I\in \D}|I|^{-2}|\langle f, \mathbbm{h}_I\rangle|^2 \mathbbm{1}_I\Big)^{1/2}\]
are the \emph{dyadic maximal function} and the \emph{dyadic square function}. Now the well-known $\textup{L}^p$ estimates for $Mf$ and $Sf$ give us the desired estimate \eqref{lp:trilinear}.

\begin{figure}[h]
\centerline{\includegraphics[scale=1]{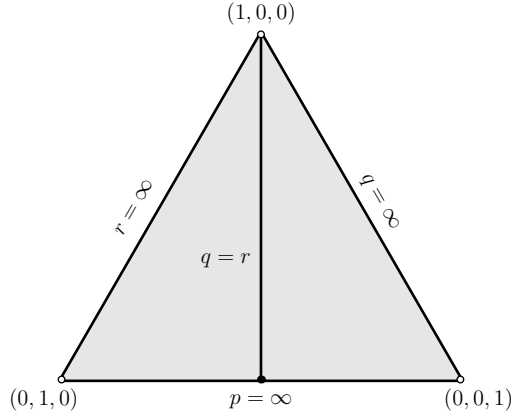}}
\caption{The Banach triangle with barycentric coordinates $(\frac 1p, \frac 1q, \frac 1r)$.}
\label{fig1}
\end{figure}

On the side $p=\infty$ of the triangle in Figure~\ref{fig1} without loss of generality we can assume that $f\equiv 1$. The sharp constant in \eqref{lp:trilinear} was found by Burkholder in \cite{Burk2} and it equals
\[\mathcal{C}_{\infty,q,r}=\max\{q-1,r-1\}.\]
On the other hand, on the sides $q=\infty$ and $r=\infty$ instead of the $\textup{L}^p$ estimates it is more natural to consider the $\textup{BMO}$ estimates, which will not be discussed in this paper. On the altitude $q=r$ of the triangle in Figure~\ref{fig1} the $\textup{L}^p$ estimates for the trilinear form \eqref{dydpar:tri} reduce to the $\textup{L}^p$ estimates for the dyadic square function, since
\[\int_\R f(Sg)^2=\Lambda_{\epsilon}(f,g,g)\]
if we take $\epsilon_I=1$ for each $I\in\D$. This implies $\|(Sg)^2\|_{\textup{L}^{q/2}(\R)} \leqslant \mathcal{C}_{p,q,q}\|g\|_{\textup{L}^q(\R)}^2$, i.e.
\[\|S\|_{\textup{L}^q(\R) \to \textup{L}^q(\R)} \leqslant \sqrt{\mathcal{C}_{p,q,q}}.\]
Actually, if the constant $\mathcal{C}_{p,q,q}$ is sharp, the last inequality turns into an equality. That sharp constant was found by Davis in \cite{Davis} and it equals
\[\mathcal{C}_{p,q,q}=(z_q^*)^{-2},\]
where $z_q^*$ is the smallest positive zero of the confluent hypergeometric function (see \cite{as:hmf}).

The special cases listed above are well-studied and even the appropriate Bellman functions are found. For $p=\infty$ one can find them in the papers by Burkholder \cite{Burk2}, Nazarov and Treil \cite{NT}, Vasyunin and Volberg \cite{VV}, Ba\~{n}uelos and Os\c{e}kowski \cite{BO}, while for $q=r$ the reader can consult the book by Os\c{e}kowski \cite{Osek}. Therefore, because of the symmetry, throughout this paper we restrict our attention to the triples of exponents $(p,q,r)$ satisfying
\begin{equation}\label{eq:exponents}
1<p,q,r<\infty,\quad q>r,\quad \frac{1}{p}+\frac{1}{q}+\frac{1}{r} = 1,
\end{equation}
which correspond to the right half of the Banach triangle depicted in Figure~\ref{fig1}.

Our goal is to give a direct proof of \eqref{lp:trilinear} using the Bellman function method. Such proofs typically give a better quantitative control and the same Bellman function can often be applied in various other settings.

First, we may assume that $f,g,h$ are non-negative, as otherwise we split them into positive and negative parts. Furthermore, we observe that it turns out to be more practical to replace the right hand side $\|f\|_{\textup{L}^p(\R)}\|g\|_{\textup{L}^q(\R)}\|h\|_{\textup{L}^r(\R)}$ by the quantity
\[\frac 1p \|f\|_{\textup{L}^p(\R)}^p +\frac 1q \|g\|_{\textup{L}^q(\R)}^q+\frac{1}{r}\|h\|_{\textup{L}^r(\R)}^r,\]
which is larger by Young's inequality, but the newly obtained inequality is actually equivalent to \eqref{lp:trilinear}, because of the homogeneity of the left hand side.

Since $|\epsilon_I |\leqslant 1$ are arbitrary, it is enough to prove a non-homogeneous estimate
\begin{equation*}
\sum_{I\in\mathcal{D}}|I|[f]_I\Big|\frac{[g]_{I\n}-[g]_{I\m}}{2}\Big |\Big| \frac{[h]_{I\n}-[h]_{I\m}}{2}\Big |\ \leqslant \mathcal{C}_{p,q,r} \Big(\frac 1p \|f\|_{\textup{L}^p(\R)}^p +\frac 1q \|g\|_{\textup{L}^q(\R)}^q+\frac{1}{r}\|h\|_{\textup{L}^r(\R)}^r\Big).
\end{equation*}
If we want to recover the original estimate \eqref{lp:trilinear}, we just have to homogenize the above inequality and use the assumed bound on $\epsilon_I$.

For an arbitrary dyadic interval $I$ we define a scale-invariant expression
\[\Phi_I(f,g,h):=\frac{1}{|I|}\sum_{\substack{J\in\D\\J \subseteq I}}|J|[f]_J\frac{|[g]_{J\n}-[g]_{J\m}|}{2}\frac{|[h]_{J\n}-[h]_{J\m}|}{2}, \]
so that we can normalize the desired estimate and rewrite it as
\begin{equation}
\Phi_I(f,g,h) \leqslant \mathcal{C}_{p,q,r}\Big( \frac 1p [f^p]_I+\frac 1q [g^q]_I+\frac 1r [h^r]_I\Big).
\label{lp:normalized}
\end{equation}
This is easily seen multiplying \eqref{lp:normalized} by $|I|$ and letting $I$ exhaust the positive and the negative half-axis.
Splitting $\sum_{J\subseteq I}$ into $\sum_{J\subseteq I\n}$, $\sum_{J\subseteq I\m}$, and $J=I$ gives us the following scaling identity:
\begin{equation}
\Phi_I(f,g,h)=\frac{1}{2}\Phi_{I\n}(f,g,h)+\frac{1}{2}\Phi_{I\m}(f,g,h)+[f]_I\frac{|[g]_{I\n}-[g]_{I\m}|}{2}\frac{|[h]_{I\n}-[h]_{I\m}|}{2}.
\label{eq:scaling}
\end{equation}
We can define the \emph{abstract Bellman function}
\[\mathbb{B}(u,v,w,U,V,W):=\sup_{f,g,h} \Phi_I(f,g,h),\]
where the supremum is taken over all non-negative functions $f,g,h$ such that
\[ [f]_I=u,  \quad [g]_I=v, \quad [h]_I=w, \quad [f^p]_I=U,  \quad [g^q]_I=V, \quad [h^r]_I=W.\]
Note that the above supremum does not depend on the choice of the ``base'' interval $I$.

Now we list some properties of that function.
\begin{enumerate}[label=($\mathcal{B}$\arabic*),ref=$\mathcal{B}$\arabic*]
\item \label{eq:b1}
\emph{Domain:} The function $\mathbb{B}$ is defined on the set
\[\mathbb{D}:=\{(u,v,w,U,V,W)\in [0,\infty)^6: u^p\leqslant U, v^q\leqslant V, w^r\leqslant W\}.\]
The upper bounds simply follow from Jensen's inequality.
\item \label{eq:b2}
\emph{Range:}
\[0\leqslant \mathbb{B}(u,v,w,U,V,W)  \leqslant \mathcal{C}_{p,q,r}\Big(\frac 1p U+\frac 1q V+\frac 1r W \Big), \]
where on the right hand side we assume that the estimate \eqref{lp:normalized} holds.
\item \label{eq:b3}
\emph{The main inequality:}
\begin{equation*}
\mathbb{B}(\mathbf{x})\geqslant \frac{1}{2}\mathbb{B}(\mathbf{x}_1)+\frac{1}{2}\mathbb{B}(\mathbf{x}_2) +u\frac{|v_1-v_2|}{2}\frac{|w_1-w_2|}{2},
\end{equation*}
whenever the six-tuples $\mathbf{x}=(u,v,w,U,V,W)$ and $\mathbf{x}_i=(u_i,v_i,w_i,U_i,V_i,W_i)$, $i=1,2$, belong to the domain and satisfy $\mathbf{x}=\frac 12 \mathbf{x}_1+\frac 12 \mathbf{x}_2$. This can be easily seen by taking the supremum in the scaling identity \eqref{eq:scaling} over all non-negative functions $f,g,h$ such that $[f]_{I\n}=u_1$, $[f^p]_{I\n}=U_1$, etc.
\end{enumerate}

Conversely, suppose that we have already found a function $\B$ with properties \eqref{eq:b1}--\eqref{eq:b3}. We will show how its existence implies the estimate $\eqref{lp:trilinear}$. Applying \eqref{eq:b3} $n$ times with a fixed choice of the functions $f,g,h \geqslant 0$ and a fixed base interval $I$ gives us
\begin{align*}
|I|\,\B\big([f]_I,[g]_I,[h]_I,[f^p]_I,[g^q]_I,[h^r]_I\big) & \geqslant \sum_{\substack{J\subseteq I\\|J|=2^{-n}|I|}} |J| \B \big([f]_J,[g]_J,[h]_J,[f^p]_J,[g^q]_J,[h^r]_J\big) \\
& + \sum_{\substack{J\subseteq I\\|J|>2^{-n}|I|}}  |J| [f]_J \ \frac{\big|[g]_{J\n}-[g]_{J\m}\big|}{2}\ \frac{\big|[h]_{J\n}-[h]_{J\m}\big|}{2}.
\end{align*}
Since by \eqref{eq:b2} the first sum is non-negative and
\[ \B\big([f]_I,[g]_I,[h]_I,[f^p]_I,[g^q]_I,[h^r]_I\big) \leqslant \mathcal{C}_{p,q,r}\Big( \frac 1p [f^p]_I+\frac 1q [g^q]_I+\frac 1r [h^r]_I\Big), \]
letting $n\to\infty$ leads us to the estimate \eqref{lp:normalized} and then in turn also to \eqref{lp:trilinear}.

It will be convenient to find a function $\B$ that also satisfies the following condition:
\begin{equation} \label{eq:b4}
\B(\mathbf{x})+(d\B)(\mathbf{x})(\mathbf{x}_1-\mathbf{x}) \geqslant \B(\mathbf{x}_1)+\frac{2}{3}u|v_1-v||w_1-w|,
\tag{$\B4$}
\end{equation}
whenever the six-tuples $\mathbf{x}=(u,v,w,U,V,W)$ and $\mathbf{x}_1=(u_1,v_1,w_1,U_1,V_1,W_1)$ belong to the domain \eqref{eq:b1}. Here $d\B$ denotes the differential of $\B$, which is a linear form, and we consider it at the point $\mathbf{x}$ and apply it to the vector $\mathbf{x}_1-\mathbf{x}$. Condition \eqref{eq:b4} is required by an application considered in Subsection \ref{subsection31}.

Now we want to find an explicit formula for one possible function $\B$. We define the function $\B\colon\mathbb{D}\to\R$ as
\begin{equation}
\B(u,v,w,U,V,W):=\mathcal{C}_{p,q,r}\Big(\frac 1p U+\frac 1q V+\frac 1r W \Big )-\A(u,v,w),
\label{eq:bellman}
\end{equation}
where $\A\colon[0,\infty)^3\to\R$ is given by
\begin{align*}
&\A(u,v,w): = \\
&\left\{
\begin{array}{ll}
      Au^p+Bv^q+Cw^r; & \, u^p\leqslant w^r \leqslant v^q,\\
      \frac{A(p-1)-C}{p-1}u^p+Bv^q+\frac{Cp}{p-1}uw^{r-\frac{r}{p}}; & \, w^r\leqslant u^p \leqslant v^q,\\
      \frac{A(p-1)-(B+C)}{p-1}u^p+\frac{Bp}{p-1}uv^{q-\frac{q}{p}}+\frac{Cp}{p-1}uw^{r-\frac{r}{p}}; & \, w^r\leqslant v^q\leqslant u^p,\\
	 \frac{A(p-1)-(B+C)}{p-1}u^p+ \frac{Bq}{2}uv^2w^{1-\frac rq}+\frac{2Cpr-Bp(q-r)}{2r(p-1)}uw^{r-\frac{r}{p}}; & \, v^q\leqslant w^r\leqslant u^p, \\
	 \frac{2Ar(p-1)-B(q+r)}{2r(p-1)}u^p+\frac{Bq^2}{2p(q-2)}u^{p-\frac{2p}{q}}v^2+\frac{Bq(q-r)}{2r(q-2)}v^2w^{r-\frac{2r}{q}}+\frac{2Cr-B(q-r)}{2r}w^r; & \, v^q\leqslant u^p\leqslant w^r, \\
	 Au^p+\frac{Bq}{p(q-2)}v^q+\frac{Bq(q-r)}{2r(q-2)}v^2w^{r-\frac{2r}{q}}+\frac{2Cr-B(q-r)}{2r}w^r; & \, u^p\leqslant v^q \leqslant w^r.
\end{array}
\right.
\end{align*}
The coefficients $A,B,C>0$ will be appropriately chosen depending only on the exponents $p,q,r$ and then one will be able to take $\mathcal{C}_{p,q,r} = \max\{Ap,Bq,Cr\}$.
We see that the function $\A$ has a similar form to the one constructed by Nazarov and Treil \cite{NT}, which can in our notation be written as
\begin{equation*}
\mathcal{NT}(v,w)=A(v^q+w^r)+B\left\{
\begin{array}{ll}
\frac 2q v^q+\big(\frac{2}{r}-1\big)w^r; & \quad v^q \geqslant w^r, \\
v^2w^{2-r}; & \quad v^q \leqslant w^r.
\end{array}
\right.
\end{equation*}
It corresponds to the endpoint case $p=\infty$, $1<r<2<q<\infty$. Instead of one critical curve $v^q=w^r$ for $\mathcal{NT}$, we have three critical surfaces:
\begin{equation}
u^p=v^q, \quad u^p=w^r, \quad v^q=w^r.
\label{hyper}
\end{equation}

Finally, we are ready to state our main result.

\begin{theorem}\label{thm:bellman}
For the exponents $p,q,r$ satisfying \eqref{eq:exponents} it is possible to \linebreak choose the coefficients $A,B,C$ such that the function $\B$ defined by \eqref{eq:bellman} is of class $\textup{C}^1$ on the whole domain $\mathbb{D}$ and satisfies the conditions \eqref{eq:b2} (with $\mathcal{C}_{p,q,r} = \max\{Ap,Bq,Cr\}$), \eqref{eq:b3}, and \eqref{eq:b4}. One possible choice of the coefficients is
\[A=\frac{88q^4r}{(p-1)(r-1)(q-r)}, \qquad B=1, \qquad \text{and} \qquad
C=\frac{11q^3r}{(r-1)(q-r)},\]
which yields
\[\mathcal{C}_{p,q,r}=\frac{88pq^4r}{(p-1)(r-1)(q-r)}.\]
\end{theorem}

The claim that $\B$ is of class $\textup{C}^1$ on $\mathbb{D}$ should be understood in the sense that the function $\A$ is continuous on $[0,\infty)^3$, $\A$ is continuously differentiable on $(0,\infty)^3$, and the partial derivatives of $\A$ can be continuously extended to $[0,\infty)^3$. At a boundary point the differential $d\B$ in \eqref{eq:b4} is interpreted as the linear form whose coefficients are the aforementioned continuous extensions of partial derivatives to that point.

The motivation behind finding the explicit Bellman function (instead of just using the abstract one) is that in some contexts the explicit formula could be useful. For example, Carbonaro and Dragi\v{c}evi\'{c} in \cite{CD1} and \cite{CD2} made use of the fact that the explicit Bellman function $\mathcal{NT}$ involves powers. Another source of motivation is that we would also like to find a direct proof (without stopping time arguments) of the estimates for the ``twisted'' paraproduct considered by one of the authors in \cite{vk:tp} or the ``twisted'' quadrilinear form considered by Durcik in \cite{pd:L4} and \cite{pd:Lp}. This could also extend the range of exponents for a non-adapted stochastic integral considered by the authors in \cite{ks:mart} or for the norm-variation of ergodic averages with respect to two commuting transformations \cite{DKST}. So far we can only say that the Bellman function that has to be constructed for any of the mentioned problems should necessarily encode some structure of the function from Theorem~\ref{thm:bellman}, as dyadic paraproducts are the simplest and prototypical multilinear multipliers.

The Bellman function that we construct certainly does not give the best possible constants $\mathcal{C}_{p,q,r}$ in \eqref{lp:trilinear}. Indeed, the sharp constant for any triple of exponents from the generic range \eqref{eq:exponents} has not yet been determined to the best of our knowledge. Search for the abstract Bellman function $\mathbb{B}$ would lead us to the equations
\begin{equation}
\det\begin{bmatrix}
\partial_u^2\mathbb{B} & \partial_u\partial_v\mathbb{B} & \partial_u\partial_w\mathbb{B} & \partial_u\partial_U\mathbb{B} & \partial_u\partial_V\mathbb{B} & \partial_u\partial_W\mathbb{B} \\
\partial_u\partial_v\mathbb{B} & \partial_v^2\mathbb{B} & \partial_v\partial_w\mathbb{B}\pm u & \partial_v\partial_U\mathbb{B} & \partial_v\partial_V\mathbb{B} & \partial_v\partial_W\mathbb{B} \\
\partial_u\partial_w\mathbb{B} & \partial_v\partial_w\mathbb{B}\pm u & \partial_w^2\mathbb{B} & \partial_w\partial_U\mathbb{B} & \partial_w\partial_V\mathbb{B} & \partial_w\partial_W\mathbb{B} \\
\partial_u\partial_U\mathbb{B} & \partial_v\partial_U\mathbb{B} & \partial_w\partial_U\mathbb{B} & \partial_U^2\mathbb{B} & \partial_U\partial_V\mathbb{B} & \partial_U\partial_W\mathbb{B} \\
\partial_u\partial_V\mathbb{B} & \partial_v\partial_V\mathbb{B} & \partial_w\partial_V\mathbb{B} & \partial_U\partial_V\mathbb{B} & \partial_V^2\mathbb{B} & \partial_V\partial_W\mathbb{B} \\
\partial_u\partial_W\mathbb{B} & \partial_v\partial_W\mathbb{B} & \partial_w\partial_W\mathbb{B} & \partial_U\partial_W\mathbb{B} & \partial_V\partial_W\mathbb{B} & \partial_W^2\mathbb{B}
\end{bmatrix}=0.
\label{PDE1}
\end{equation}
One way of simplifying \eqref{PDE1} is to consider the non-homogeneous function $\B$ of the form \eqref{eq:bellman}. Function $\mathcal{B}$ is now a supersolution of the equation for the true Bellman function $\mathbb{B}$, but a function of that form can still yield the optimal (unknown) constant. This way \eqref{PDE1} reduces to  
\begin{equation}
\det \mathbb{A}_{\pm}=0,
\label{PDE2}
\end{equation}
where $\mathbb{A}_{\pm}$ are the matrices defined with \eqref{matrix} below. Alternatively, one can use the homogeneities of $\mathbb{B}$ to reduce the dimension in \eqref{PDE1}. Equations like \eqref{PDE2} can sometimes be turned into the Monge-Amp\` ere equation by an appropriate change of variables, which does not seem to be the case here. At the moment, we do not know how to solve \eqref{PDE2}, so we impose slightly weaker conditions on our function $\B$ that result in a constant $\mathcal{C}_{p,q,r}$ which is not optimal. It would be interesting to find a Bellman function $\mathcal{B}$ that yields the optimal constant, or perhaps even the exact abstract Bellman function $\mathbb{B}$.  Let us remark once again that this was achieved by Ba\~{n}uelos and Os\c{e}kowski \cite{BO} in the endpoint case $p=\infty$, $f\equiv 1$.

We have organized the remainder of the paper as follows. In the next section we present the proof of Theorem~\ref{thm:bellman}. In Section~\ref{section3} we apply Theorem~\ref{thm:bellman} to reprove the well-known $\textup{L}^p$ estimates for martingale paraproducts and the heat flow paraproducts.


\section{Proof of Theorem~\ref{thm:bellman}}
The continuity of $\A$ on $[0,\infty)^3$ is obvious. Indeed, observe that all exponents appearing in the definition of $\mathcal{A}$ are positive. Thus, $\mathcal{A}$ is clearly well-defined and continuous on each of the six closed regions determined by the inequalities for $u,v,w$ and it is straightforward to verify that the six formulas are compatible on the common boundaries.

To see that $\A$ is continuously differentiable on the open octant $(0,\infty)^3$ we just calculate the first order partial derivatives in the interior of each of the previously mentioned regions. The formula for each of these derivatives inside any of the regions continuously extends to the whole open octant. Moreover, these formulas coincide on the boundaries of each two adjacent regions, so we can deduce that $\A$ really is of class $\textup{C}^1$ on $(0,\infty)^3$. For instance, when we calculate the partial derivative of $\A$ with respect to $u$ on the two adjacent open regions $v^q<u^p<w^r$ and $u^p<v^q<w^r$, we get
\[\frac{\partial \A}{\partial u}(u,v,w)=\frac{2Arp(p-1)-Bp(q+r)}{2r(p-1)}u^{p-1}+\frac{Bq}{2}u^{\frac pr-\frac pq}v^2 \quad \text{and} \quad \frac{\partial \A}{\partial u}(u,v,w)=Apu^{p-1}.\]
The common boundary of those two regions is a subset of $u^p=v^q$, where both formulas give $Apu^{p-1}$, i.e.\@ the two formulas for the partial derivative coincide on that boundary. All other cases are treated in the same manner.

Also, it is easy to see that the partial derivatives have limits at each point of the boundary of $[0,\infty)^3$ and hence they can be continuously extended to $[0,\infty)^3$. For example, if $0<v^q\leqslant w^r \leqslant u^p$, then the partial derivative of $\A$ with respect to $w$ equals
\[\frac{\partial \A}{\partial w}(u,v,w)=\frac{B(q-r)}{2}uv^2w^{-\frac rq}+\frac{2Cr-B(q-r)}{2}uw^{\frac{r}{q}}.\]
Obviously, the only problematic points are the ones on the part of the boundary lying on the plane $w=0$, but since $\frac{v^q}{w^r}\leqslant 1$, we have
\begin{equation*}
\displaystyle \lim_{w \to 0+} \frac{\partial \A}{\partial w}(u,v,w) = \lim_{w \to 0+} \left(\frac{B(q-r)}{2}uw^\frac{r}{q}\Big(\frac{v^q}{w^r}\Big)^{\frac 2q}+\frac{2Cr-B(q-r)}{2}uw^{\frac{r}{q}}\right)=0.
\end{equation*}
We can show the existence of the other limits in a similar way.

The estimate \eqref{eq:b2} follows directly from the definitions of the functions $\A$ and $\B$, since
\begin{equation} \label{eq:a2}
0\leqslant \A(u,v,w)\leqslant Au^p+Bv^q+Cw^r
\tag{$\A2$}
\end{equation}
as long as $A,B,C\geqslant 0$.
This is easily seen using Young's inequality. For instance, if $w^r\leqslant v^q\leqslant u^p$, then we have
\begin{align*}
\A(u,v,w)&=\frac{A(p-1)-(B+C)}{p-1}u^p+\frac{Bp}{p-1}uv^{q-\frac qp}+\frac{Cp}{p-1}uw^{r-\frac rp} \\
&\leqslant \frac{A(p-1)-(B+C)}{p-1}u^p+\frac{Bp}{p-1}\Big(\frac 1p u^p+\frac{p-1}{p}v^q\Big)+\frac{Cp}{p-1}\Big(\frac{1}{p}u^p+\frac{p-1}{p}w^r\Big).
\end{align*}
Other cases follow analogously.
Non-negativity of $\B$ on $\mathbb{D}$ is guaranteed if $\mathcal{C}_{p,q,r}\geqslant Ap,Bq,Cr$.

Observe that \eqref{eq:b3} is equivalent to
\begin{equation} \label{eq:a3}
\frac{1}{2}\A(u_1,v_1,w_1)+\frac{1}{2}\A(u_2,v_2,w_2)-\A(u,v,w) \geqslant u\frac{|v_1-v_2|}{2}\frac{|w_1-w_2|}{2},
\tag{$\A3$}
\end{equation}
where $(u,v,w)$, $(u_1,v_1,w_1)$, and $(u_2,v_2,w_2)$ are in $[0,\infty)^3$ and such that
\begin{equation}
(u,v,w)=\frac{1}{2}(u_1,v_1,w_1)+\frac{1}{2}(u_2,v_2,w_2),
\label{a3cond}
\end{equation}
while \eqref{eq:b4} is equivalent to
\begin{equation} \label{eq:a4}
\A(u_1,v_1,w_1) \geqslant \A(u,v,w)+(d\A)(u,v,w)(u_1-u,v_1-v,w_1-w)+\frac 23 u|v_1-v||w_1-w|,
\tag{$\A 4$}
\end{equation}
where $(u,v,w)$ and $(u_1,v_1,w_1)$ are in $[0,\infty)^3$. Instead of proving \eqref{eq:a3} and \eqref{eq:a4} directly, we will reduce them conveniently to an inequality for quadratic forms.

Let $(u,v,w)\in (0,\infty)^3$ be a point that does not lie on any of the three critical surfaces \eqref{hyper}. This means that $\A$ is of class $\textup{C}^2$ on an open ball around that point. If we take $(u_1,v_1,w_1), (u_2,v_2,w_2)$ from that open ball such that \eqref{a3cond} holds,
then substituting $u=\frac{u_1+u_2}{2}, \Delta u=\frac{u_1-u_2}{2}$, etc., and adding Taylor's formulas at $(u,v,w)$ for $\A(u\pm\Delta u,v\pm\Delta v,w\pm\Delta w)$ gives us the infinitesimal version of \eqref{eq:a3}:
\begin{equation}
(d^2 \A)(u,v,w)(\Delta u, \Delta v, \Delta w) \geqslant 2u|\Delta v||\Delta w|.
\label{eq:a3inf}
\tag{$\A3'$}
\end{equation}
Here $d^2\A$ denotes the second differential of $\A$ as a quadratic form, which we consider at the point $(u,v,w)$ and apply to the vector $(\Delta u, \Delta v, \Delta w)$. Notice that \eqref{eq:a3inf} does not hold on the whole domain of the function $\A$, which is $[0,\infty)^3$, but it does hold on the interior of each of the six regions into which the three surfaces divide $(0,\infty)^3$.

Conversely, \eqref{eq:a3inf} implies \eqref{eq:a3}, i.e.\@ the two inequalities are equivalent for continuously differentiable functions, which is enabled by the convexity of the domain. To show the converse, first take a point $(u,v,w) \in (0,\infty)^3$ and a vector $(\Delta u, \Delta v, \Delta w) \in \R^3$ such that also $(u\pm\Delta u, v\pm\Delta v, w\pm\Delta w)\in (0,\infty)^3$. Now define the function $\alpha\colon[-1,1]\to \R$ as
\begin{equation}
\alpha(t):=\A(u+t\Delta u, v+t\Delta v, w+t\Delta w).
\label{alfa}
\end{equation}
This function is continuously differentiable on $[-1,1]$ since $\A$ is of class $\textup{C}^1$ on $(0,\infty)^3$. Also, $\alpha$ is piecewise $\textup{C}^2$ on $[-1,1]$. This follows from the facts that $\A$ is of class $\textup{C}^2$ on $(0,\infty)^3$ outside the surfaces \eqref{hyper}, it has bounded second derivatives away from the coordinate planes $u=0$, $v=0$, and $w=0$, and the segment
$\{(u+t\Delta u, v+t\Delta v, w+t\Delta w): t\in[-1,1]\}$
intersects the three critical surfaces at finitely many points. If we denote those points $t_1<t_2<\dots <t_n$, then using the integration by parts and the fundamental theorem of calculus (both in the versions for absolutely continuous functions; see \cite{Cohn}) gives us
\[\displaystyle\frac{1}{2}\alpha(1)+\frac{1}{2}\alpha(-1)-\alpha(0)=\frac{1}{2}\int_{[-1,1]\setminus\{t_1,\dots,t_n\}}(1-|t|)\alpha''(t)dt.\]
From the above identity we deduce
\begin{gather*}
\frac 12 \A(u+\Delta u,v+\Delta v,w+\Delta w)+\frac 12\A(u-\Delta u,v-\Delta v,w-\Delta w)-\A(u,v,w) \\
=\frac{1}{2}\int_{[-1,1]\setminus\{t_1,\dots,t_n\}}(1-|t|)(d^2 \A)(u+t\Delta u, v+t\Delta v, w+t\Delta w)(\Delta u, \Delta v, \Delta w)dt.
\end{gather*}
Finally, \eqref{eq:a3inf} implies that the last expression is at least
\begin{equation*}
\frac{1}{2}\int_{[-1,1]\setminus\{t_1,\dots,t_n\}}(1-|t|)2(u+t\Delta u)|\Delta v||\Delta w|dt=u|\Delta v||\Delta w|,
\end{equation*}
which gives exactly \eqref{eq:a3}.

Moreover, \eqref{eq:a3inf} also implies \eqref{eq:a4}. In order to verify this we also take $(u,v,w)\in (0,\infty)^3$ and $(\Delta u, \Delta v, \Delta w)\in \R^3$ such that $(u+\Delta u, v+\Delta v, w+\Delta w)\in (0,\infty)^3$. We define $\alpha\colon[0,1]\to\R$ again by the formula \eqref{alfa}. Integration by parts, the fundamental theorem of calculus, and \eqref{eq:a3inf} this time give
\[\alpha(1)=\alpha(0)+\alpha'(0)+\int_{[0,1]\setminus\{t_1,\dots,t_n\}}(1-t)\alpha''(t)dt,\]
and therefore,
\begin{gather*}
\A(u+\Delta u,v+\Delta v,w+\Delta w)=\A(u,v,w)+(d\A)(u,v,w)(\Delta u, \Delta v, \Delta) \\
+\int_{[0,1]\setminus\{t_1,\dots,t_n\}}(1-t)(d^2 \A)(u+t\Delta u, v+t\Delta v, w+t\Delta w)(\Delta u, \Delta v, \Delta w)dt \\
\geqslant \A(u,v,w)+(d\A)(u,v,w)(\Delta u, \Delta v, \Delta w)+\int_{[0,1]\setminus\{t_1,\dots,t_n\}}(1-t)2(u+t\Delta u)|\Delta v||\Delta w|dt.
\end{gather*}
Since $u+t\Delta u= (1-t)u+t(u+\Delta u)\geqslant (1-t)u$, integrating in $t$ in the above inequality yields
\begin{equation*}
\A(u+\Delta u,v+\Delta v,w+\Delta w)\geqslant \A(u,v,w)+(d\A)(u,v,w)(\Delta u, \Delta v, \Delta w)+\frac 23 u|\Delta v||\Delta w|,
\end{equation*}
which is exactly \eqref{eq:a4}.

This way we proved that \eqref{eq:a3inf} implies \eqref{eq:a3} and \eqref{eq:a4}, but only on $(0,\infty)^3$. To see that these two also hold on $[0,\infty)^3$ we just have to extend the obtained inequalities by the continuity of $\A$ and $d\A$. We have commented in the introduction how we interpret $d\mathcal{A}$ at the boundary of the domain.

Now we are left with proving \eqref{eq:a3inf}, which is equivalent to showing that the two matrices
\begin{equation}
\mathbb{A}_{\pm}=\begin{bmatrix}
\partial_u^2\A & \partial_u\partial_v\A & \partial_u\partial_w\A\\
\partial_u\partial_v\A & \partial_v^2\A & \partial_v\partial_w\A\pm u\\
\partial_u\partial_w\A & \partial_v\partial_w\A\pm u & \partial_w^2\A
\end{bmatrix}
\label{matrix}
\end{equation}
are positive semi-definite on each of the six open regions into which the surfaces \eqref{hyper} split $(0,\infty)^3$. To do so we will use Sylvester's criterion and verify that all three principal minors are positive. More precisely, we will prove that the constants $A,B,C$ can be chosen so that this is fulfilled.

We can simplify the calculations a bit by substituting $t=\frac{v^q}{u^p}$, $s=\frac{w^r}{u^p}$ and noting that
\begin{equation}
\A(u,v,w)=u^p \gamma(t,s),
\label{eq:subst}
\end{equation}
where $\gamma\colon(0,\infty)^2\to\R$ is given by
\begin{equation*}
\gamma(t,s) = \left\{
\begin{array}{ll}
      A+Bt+Cs; & \quad 1\leqslant s \leqslant t,\\
      \frac{A(p-1)-C}{p-1}+Bt+ \frac{Cp}{p-1}s^{1-\frac{1}{p}}; & \quad s\leqslant 1 \leqslant t,\\
      \frac{A(p-1)-(B+C)}{p-1}+ \frac{Bp}{p-1}t^{1-\frac{1}{p}}+\frac{Cp}{p-1}s^{1-\frac{1}{p}}; & \quad s\leqslant t\leqslant 1,\\
	 \frac{A(p-1)-(B+C)}{p-1}+ \frac{Bq}{2}t^{\frac 2q}s^{\frac 1r-\frac 1q}+ \frac{2Cpr-Bp(q-r)}{2r(p-1)}s^{1-\frac{1}{p}}; & \quad t\leqslant s\leqslant 1, \\
	 \frac{2Ar(p-1)-B(q+r)}{2r(p-1)}+\frac{Bq^2}{2p(q-2)}t^{\frac{2}{q}}+\frac{Bq(q-r)}{2r(q-2)}t^{\frac{2}{q}}s^{1-\frac 2q}+\frac{2Cr-B(q-r)}{2r}s; & \quad t\leqslant 1\leqslant s, \\
	 A+\frac{Bq}{p(q-2)}t+\frac{Bq(q-r)}{2r(q-2)}t^{\frac{2}{q}}s^{1-\frac 2q}+\frac{2Cr-B(q-r)}{2r}s; & \quad 1\leqslant t \leqslant s.
\end{array}
\right.
\end{equation*}

After plugging \eqref{eq:subst} into \eqref{matrix} and multiplying from both sides with the diagonal matrix $\textup{diag}(u^{1-p/2},\, u^{p/q-p/2},\, u^{p/r-p/2})$ we obtain the matrices $M=[m_{ij}]$, where
{\allowdisplaybreaks
\begin{align*}
m_{11}&=p(p-1)\gamma(t,s)-p(p-1)t\partial_t\gamma(t,s) -p(p-1)s\partial_s\gamma(t,s),\\
&+2p^2ts\partial_t\partial_s\gamma(t,s)+p^2t^2\partial_t^2\gamma(t,s) +p^2s^2\partial_s^2\gamma(t,s), \\
m_{12}&=m_{21}=-p q t^{1-\frac{1}{q}}s\partial_t\partial_s \gamma(t,s)-pq t^{2-\frac{1}{q}} \partial_t^2\gamma(t,s), \\
m_{13}&=m_{31}=-prt s^{1-\frac{1}{r}} \partial_t\partial_s \gamma(t,s)-p r s^{2-\frac{1}{r}} \partial_s^2\gamma(t,s),\\
m_{22}&=q(q-1)t^{1-\frac{2}{q}}\partial_t \gamma(t,s)+q^2 t^{2-\frac{2}{q}} \partial_t^2\gamma(t,s),\\
m_{23}&=m_{32}=q r t^{1-\frac{1}{q}} s^{1-\frac{1}{r}} \partial_t\partial_s\gamma(t,s) \pm 1,\\
m_{33}&=r(r-1)s^{1-\frac{2}{r}} \partial_s\gamma(t,s)+r^2 s^{2-\frac{2}{r}} \partial_s^2\gamma(t,s),
\end{align*}}
and the problem is reduced to verifying that these matrices are positive definite on the interior of each of the six regions determined by the inequalities for $t$ and $s$.
First, we will calculate the three principal minors of the above matrices for each region, and then we will explain why we can choose the constants $A,B,C$ such that all of them are positive.

The following expressions were calculated using Mathematica \cite{w:m9}.
\begin{itemize}
\item[Region\!] 1: $1< s < t$
\begin{itemize}
\item[] Minor $1\times 1$: $\quad \doublebox{$Ap(p-1)$}$
\item[] Minor $2 \times 2$: $\quad \doublebox{$ABp(p-1)q(q-1)t^{1-\frac{2}{q}}$}$ \vskip 2mm
\item[] Determinants (with $\pm$):
\begin{equation*}
\doublebox{$ABCp(p-1)q(q-1)r(r-1)t^{1-\frac{2}{q}} s^{1-\frac{2}{r}}$}-Ap(p-1)
\end{equation*}
\end{itemize}
\item[Region\!] 2: $s< 1 < t$
\begin{itemize}
\item[] Minor $1\times 1$: $\quad \doublebox{$p(A(p-1)-C)$}$
\item[] Minor $2 \times 2$: $\quad \doublebox{$Bp(A(p-1)-C)q(q-1)t^{1-\frac{2}{q}}$}$ \vskip 2mm
\item[] Determinants (with $\pm$):
\begin{equation*}
\doublebox{$BCp(A(p-1)-C)(q-1)r^2t^{1-\frac{2}{q}}s^{\frac{1}{q}-\frac{1}{r}}$}-BC^2q(q-1)r^2t^{1-\frac{2}{q}}s^{\frac{2}{q}}-p(A(p-1)-C)
\end{equation*}
\end{itemize}
\item[Region\!] 3: $s< t < 1$
\begin{itemize}
\item[] Minor $1\times 1$: $\quad \doublebox{$p(A(p-1)-B-C)$}$
\item[] Minor $2 \times 2$:
\begin{equation*}
\doublebox{$\displaystyle\frac{Bp(A(p-1)-B-C)q^2}{r}t^{\frac{1}{r}-\frac{1}{q}}$}-B^2q^2t^{\frac{2}{r}}
\end{equation*}
\item[] Determinants (with $\pm$):
\begin{gather*}
\doublebox{$BCp(Ap(q+r)-qr(B+C))t^{\frac{1}{r}-\frac{1}{q}}s^{\frac{1}{q}-\frac{1}{r}}$} -B^2Cqr^2t^{\frac{2}{r}}s^{\frac{1}{q}-\frac{1}{r}}-BC^2q^2rt^{\frac{1}{r}-\frac{1}{q}}s^{\frac{2}{q}}\\
-p(A(p-1)-B-C)\pm 2BCqrt^{\frac{1}{r}}s^{\frac{1}{q}}
\end{gather*}
\end{itemize}
\item[Region\!] 4: $t< s < 1$
\begin{itemize}
\item[] Minor $1\times 1$: $\quad \doublebox{$p(A(p-1)-B-C)$}$
\item[] Minor $2 \times 2$: $\quad \doublebox{$Bp(A(p-1)-B-C)qs^{\frac{1}{r}-\frac{1}{q}}$}-B^2q^2t^{\frac{2}{q}}s^{\frac 2r-\frac 2q}$ \vskip 2mm
\item[] Determinants (with $\pm$):
\begin{gather*}
\doublebox{$\displaystyle\frac{1}{2}(A(p-1)-B-C)\big(Bpr(2Cr-B(q-r))-2p\big)$} -\frac{1}{4}Bq(2Cr-B(q-r))^2s^{\frac 1q+\frac 1r}\\
\mp 2Bp(A(p-1)-B-C)(q-r)t^{\frac 1q}s^{-\frac 1q}\pm Bq(2Cr-B(q-r))t^{\frac 1q}s^{\frac 1r}\\
-\frac{1}{2}B^2p(A(p-1)-B-C)(q-r)(2q-r)t^{\frac 2q}s^{-\frac 2q}+\frac{1}{2}B^2q(2Cr-B(q-r))(q-2r)t^{\frac 2q}s^{\frac 1r-\frac 1q}\\
\pm B^2q(q-r)t^{\frac{3}{q}}s^{\frac 1r-\frac 2q}+\frac{1}{4}B^3q(q-r)(3q-r)t^{\frac 4q}s^{\frac 1r-\frac 3q}
\end{gather*}
\end{itemize}
\item[Region\!] 5: $t< 1 < s$
\begin{itemize}
\item[] Minor $1\times 1$: \[\doublebox{$\displaystyle\frac{p(2Ar(p-1)-B(q+r))}{2r}$}+\frac{Bp(q-r)}{2r}t^{\frac 2q}\]
\item[] Minor $2 \times 2$:
\begin{gather*}
\doublebox{$\displaystyle\frac{Bpq(q-r)(2Ar(p-1)-B(q+r))}{2r^2(q-2)}s^{1-\frac 2q}$}-\frac{B^2q^2(pq+q-2p)}{2p(q-2)}t^{\frac 2q}\\
+\frac{B^2pq(q-r)^2}{2r^2(q-2)}t^{\frac 2q}s^{1-\frac 2q}+\ovalbox{$\displaystyle\frac{Bq^2(2Ar(p-1)-B(q+r))}{2r(q-2)}$}
\end{gather*}
\item[] Determinants (with $\pm$):
{\allowdisplaybreaks
\begin{gather*}
\doublebox{$\displaystyle\frac{B(q-r)(p+q)(2Ar(p-1)-B(q+r))(2Cr-B(q-r))}{4r(q-2)}s^{\frac 2p}$}-\frac{Bp(q-r)}{2r}t^{\frac 2q} \\
+\ovalbox{$\displaystyle\frac{B(2Ar(p-1)-B(q+r))(2Cr-B(q-r))q^2(r-1)}{4r(q-2)}s^{1-\frac 2r}$}\\
-\frac{p(2Ar(p-1)-B(q+r))}{2r}-\frac{B^2(q-r)^2(2Ar(p-1)-B(q+r))(2pq-3p-q)}{4r(q-2)}t^{\frac 2q}s^{\frac 2p-\frac 2q} \\
\mp \frac{B^2p(q-r)^2}{r}t^{\frac 3q}s^{\frac 1p-\frac 1q} -\frac{B^2q^2(pq-2p+q)(r-1)(2Cr-B(q-r))}{4p(q-2)}t^{\frac 2q}s^{1-\frac 2r}\\
-\frac{B^3p(q-r)^3(qr-2r+q)}{4r^2(q-2)}t^{\frac 4q}s^{\frac 2p-\frac 2q}+\frac{B^2pq(q-r)^2(r-1)(2Cr-B(q-r))}{4r^2(q-2)}t^{\frac 2q}s^{\frac 2p}\\
\mp \frac{Bp(q-r)(2Ar(p-1)-B(q+r))}{r}t^{\frac 1q}s^{\frac 1p-\frac 1q}
-\frac{B^3qr(pq-2p+q)(q-r)(q-p)}{4p^2(q-2)}t^{\frac 4q}s^{\frac 2p-1}\\
+\frac{B^2q(q-r)(q-p)(2Ar(p-1)-B(q+r))}{4p(q-2)}t^{\frac 2q}s^{\frac 2p-1}
\end{gather*}}
\end{itemize}
\item[Region\!] 6: $1< t < s$
\begin{itemize}
\item[] Minor $1\times 1$: $\quad \doublebox{$Ap(p-1)$}$
\item[] Minor $2 \times 2$: \[\doublebox{$\displaystyle\frac{ABpq(p-1)(q-r)}{r(q-2)}s^{1-\frac 2q}$}+\ovalbox{$\displaystyle\frac{ABq^2(p-1)(q-1)}{q-2}t^{1-\frac{2}{q}}$}\]
\item[] Determinants (with $\pm$):
\begin{gather*}
\doublebox{$\displaystyle\frac{AB(2Cr-B(q-r))(p-1)(q-r)(p+q)}{2(q-2)}s^{\frac 2p}$}-\frac{AB^2qr(p-1)(q-1)(q-r)(p-q)}{2p(q-2)}ts^{\frac 2p-1}\\
+\ovalbox{$\displaystyle\frac{AB(2Cr-B(q-r))qr(p-1)(q-1)(p+q)}{2p(q-2)}t^{1-\frac 2q}s^{1-\frac 2r}$}-Ap(p-1) \\
-\frac{AB^2(p-1)(q-r)^2(2pq-3p-q)}{2(q-2)}t^{\frac 2q}s^{\frac 2p-\frac 2q}\mp 2ABp(p-1)(q-r)t^{\frac 1q}s^{\frac 1p-\frac 1q}
\end{gather*}
\end{itemize}
\end{itemize}

In each of the expressions there is a unique dominant term (regarding the exponents of $t$ and $s$) and it is double framed. We choose $B$ arbitrarily (say $B=1$), then take $C$ large enough (depending on $p,q,r,B$), and finally take $A$ large enough (depending on $p,q,r,B,C$). While doing so, we take care that the coefficient of the double framed term is greater than the sum of the absolute values of coefficients of the terms that are neither framed nor circled. We can do so because by taking $C$ large enough the expression multiplying $A$ in the coefficient of the dominant term can be made larger than the sum of the absolute values of the corresponding expressions in other non-circled terms that contain $A$.   Consequently, the coefficient of the dominant term grows faster than the sum of the absolute values of the coefficients in the other terms as $A$ tends to infinity. This means that we can take $A$ large enough so that the dominant term actually dominates the sum of all other non-framed and non-circled terms in each expression. Another way of phrasing the argument that sufficiently large $A$ and $C$ make six considered determinantal expressions positive is to observe that each dominant term contains the product $AC$, as opposed to any other non-circled term.

The only problematic terms that we cannot dominate with the dominant term are the circled ones, because of their uncontrollable growth in $A$. However, just by taking
\[C\geqslant \frac{B(q-r)}{2r} \quad \text{and} \quad A\geqslant \frac{B(q+r)}{2r(p-1)}\]
we make sure that all of them are non-negative, so they only contribute to the positivity of the expressions.

To explain how the values of the coefficients $A$, $B$, and $C$ in Theorem \ref{thm:bellman} were obtained, let us consider Region $4$ as a representative example. The other regions are treated similarly.

First, notice that the double framed term really is the dominant one, since $t<s<1$ implies
\[t^{\frac 3q}s^{\frac 1r-\frac 2q},\,
t^{\frac 1q}s^{-\frac 1q},\,
t^{\frac 4q}s^{\frac 1r-\frac 3q}\,
t^{\frac 2q}s^{-\frac 2q},\,
t^{\frac 1q}s^{\frac 1r},\,
s^{\frac 1q+\frac 1r},\,
t^{\frac 2q}s^{\frac 1r-\frac 1q}\,
<\, 1 = t^0 s^0.\]
We can choose $B=1$ and then take $C$ large enough such that
\[ r(2Cr-q+r) > \max\{ 28(q-r)+2, 7(q-r)(2q-r)+2 \}. \]
Clearly, $C=\frac{11q^3r}{(r-1)(q-r)}$ satisfies the above condition. This way the expression multiplying $A$ in the coefficient of the dominant term is seven times larger than the expressions multiplying $A$ in the coefficients of the two non-framed terms that contain $A$. Now we just have to take $A$ large enough such that 
\[(A(p-1)-C-1)\big(pr(2Cr-q+r)-2p\big)\]
is at least
\begin{equation*}
\max\Bigg\{\frac{7}{2}q(2Cr-q+r)^2,\, 14q(2Cr-q+r),\, 7q(2Cr-q+r)|q-2r|,\, 14q(q-r), \, \frac{7}{2}q(q-r)(3q-r)\Bigg\}.
\end{equation*}
It is easy to see that $A=\frac{88q^4r}{(p-1)(r-1)(q-r)}$ is one possible choice. Now the dominant term is more than seven times larger than the absolute value of any other term, which means that the dominant term dominates the sum of all other terms.

This way we accomplish the positivity of each of the expressions, which is exactly what we needed and the proof of \eqref{eq:a3inf} is completed. This also finishes the proof of Theorem~\ref{thm:bellman}.

In the next section, it will sometimes be more convenient to use the infinitesimal version of \eqref{eq:b3}:
\begin{equation} \label{eq:b3inf}
-(d^2\mathcal{B})(u,v,w,U,V,W)(\Delta u,\Delta v,\Delta w,\Delta U,\Delta V,\Delta W)
\geqslant 2u |\triangle v| |\triangle w|.
\tag{$\B 3'$}
\end{equation}
Again, \eqref{eq:b3inf} holds only for points $(u,v,w,U,V,W)$ at which the second differential of $\B$ is well-defined, i.e.\@ for the points such that $(u,v,w)$ does not lie on any of the three critical surfaces. The equivalence of \eqref{eq:b3inf} and \eqref{eq:b3} follows from the equivalence of \eqref{eq:a3inf} and \eqref{eq:a3}.

\section{Applications}
\label{section3}
Here we present several applications of the existence of the Bellman function from Theorem~\ref{thm:bellman}. We need to emphasize that the following problems are quite classical and can be solved using more standard tools. We only provide quite straightforward solutions based on Theorem~\ref{thm:bellman}. Moreover, only the existence of the Bellman function with properties \eqref{eq:b1}--\eqref{eq:b3} is needed, even though \eqref{eq:b4} is quite convenient in Subsection~\ref{subsection31}. This existence can also follow if boundedness of the dyadic paraproduct is established in some other way, as commented in the introduction. However, our goal is to illustrate how several classical problems become methodologically simple once we explicitly construct the function as in Theorem~\ref{thm:bellman}.

For two non-negative quantities $A$ and $B$ we will write $A \lesssim_P B$ if there exists a finite constant $C_P \geqslant 0$ depending on a set of parameters $P$ such that $A \leqslant C_P B$.

\subsection{Discrete-time martingales}
\label{subsection31}
Let us consider two martingales $X=(X_n)_{n=0}^\infty$ and $Y=(Y_n)_{n=0}^\infty$ with respect to the same filtration $(\mathcal{F}_n)_{n=0}^\infty$. Their \emph{paraproduct} is a stochastic process $\big((X\cdot Y)_n\big)_{n=0}^\infty$ defined as
\begin{align}
(X\cdot Y)_0&:=0, \nonumber \\
(X\cdot Y)_n&:=\sum_{k=1}^n X_{k-1}(Y_k-Y_{k-1}) \quad \text{for $n\geqslant 1$}.
\label{Dis:paraprod}
\end{align}
This process can be regarded as a particular case of Burkholder's martingale transform \cite{Burk1} of the martingale $Y$ with respect to the shifted adapted process $X$. We have also imposed the martingale property on $X$, since we want to treat $X$ and $Y$ symmetrically and since this is required by the existence of the $\textup{L}^p$ estimates in the interior of the Banach triangle in Figure~\ref{fig1}.
We want to prove that for the exponents $p,q,r$ satisfying \eqref{eq:exponents} the estimate
\begin{equation}
\|(X \cdot Y)_n\|_{\textup{L}^{r'}} \lesssim_{p,q,r} \|X_n\|_{\textup{L}^p}\|Y_n\|_{\textup{L}^q}
\label{dis:est}
\end{equation}
holds uniformly in the positive integer $n$, where $r'$ is the conjugate exponent of $r$. Instead of proving \eqref{dis:est} directly, we will rather show the estimate for the dualized form, i.e.\@ that for an arbitrary random variable $Z\in \textup{L}^r$ the inequality
\begin{equation}
|\E\big((X\cdot Y)_n Z\big)| \lesssim_{p,q,r} \|X_n\|_{\textup{L}^p}\|Y_n\|_{\textup{L}^q}\|Z\|_{\textup{L}^r}
\label{dis:est2}
\end{equation}
holds. This inequality is trivial unless all norms on the right hand side are finite.

Let us introduce the third martingale $(Z_n)_{n=0}^\infty$ with $Z_n:=\E(Z|\F_n)$. By splitting
\[Z=Z_{k-1}+(Z_k-Z_{k-1})+(Z-Z_k),\]
we can write
\begin{gather*}
\E\big((X\cdot Y)_n Z\big) =\sum_{k=1}^n \E\big(X_{k-1}(Y_k-Y_{k-1})Z\big)=\sum_{k=1}^n \E\big(X_{k-1}Z_{k-1}\E(Y_k-Y_{k-1}|\F_{k-1})\big) \nonumber\\
+\sum_{k=1}^n \E\big(X_{k-1}(Y_k-Y_{k-1})(Z_k-Z_{k-1})\big)+\sum_{k=1}^n \E\big(X_{k-1}(Y_k-Y_{k-1})\E(Z-Z_k|\F_k)\big) \nonumber \\
=\sum_{k=1}^n \E\big(X_{k-1}(Y_k-Y_{k-1})(Z_k-Z_{k-1})\big).
\end{gather*}
Here the third equality follows from the martingale property, since $\E(Y_k-Y_{k-1}|\F_{k-1})=0$ and $\E(Z-Z_k|\F_k)=0$. The estimate \eqref{dis:est2} is well-known and its proof uses the Cauchy-Schwarz, H\"{o}lder's, Doob's and the Burkholder-Gundy inequalities. Again, we will give a more direct proof using the Bellman function \eqref{eq:bellman}.

It is enough to consider the times $k=0,1,\dots,n$, but we need to show the estimate that is uniform in $n$. We can assume that $X_k,Y_k,Z_k\geqslant 0$ for $0\leqslant k \leqslant n$, as otherwise we split the variables $X_n,Y_n,Z_n$ into positive and negative parts, and introduce three new martingales (for a fixed $n$):
\[ U_k := \mathbb{E}(X_n^p|\mathcal{F}_k),\quad V_k := \mathbb{E}(Y_n^q|\mathcal{F}_k),\quad W_k := \mathbb{E}(Z_n^r|\mathcal{F}_k).  \]
If we denote $\mathbf{X}_k = (X_k,Y_k,Z_k,U_k,V_k,W_k)$, then property \eqref{eq:b4} of the Bellman function $\B$ gives us
\begin{equation*}
\B(\mathbf{X}_{k-1})+(d\B)(\mathbf{X}_{k-1})(\mathbf{X}_k-\mathbf{X}_{k-1}) \geqslant \B(\mathbf{X}_k)+\frac{2}{3}X_{k-1}|Y_k-Y_{k-1}||Z_k-Z_{k-1}|,
\end{equation*}
from which we deduce
\begin{equation}
\B(\mathbf{X}_{k-1})\geqslant \E\big(\B(\mathbf{X}_k)\big|\F_{k-1}\big)
+\frac{2}{3}\E\big(X_{k-1}|Y_k-Y_{k-1}||Z_k-Z_{k-1}|\,\big|\F_{k-1}\big),
\label{dis:ineq}
\end{equation}
by taking the conditional expectation with respect to $\F_{k-1}$ and using the martingale property. Finally, taking the expectation of \eqref{dis:ineq}, summing over $k=1,\dots,n$, telescoping, and using \eqref{eq:b2} gives
\begin{align*}
\frac{2}{3}\sum_{k=1}^n\E\big(X_{k-1}|Y_k-Y_{k-1}||Z_k-Z_{k-1}|\big) &\leqslant \E\B(\mathbf{X}_0)-\E\B(\mathbf{X}_n)\\
\leqslant \mathcal{C}_{p,q,r} \E\Big( \frac{1}{p}U_0 + \frac{1}{q}V_0 + \frac{1}{r}W_0 \Big)&=\mathcal{C}_{p,q,r}\Big(\frac 1p  \|X_n\|_{\textup{L}^p}^p+\frac 1q \|Y_n\|_{\textup{L}^q}^q+\frac 1r \|Z_n\|_{\textup{L}^r}^r \Big ).
\end{align*}
Homogenizing the above inequality we get the desired estimate $\eqref{dis:est2}$ and hence also $\eqref{dis:est}$.

\subsection{Continuous-time martingales}
Let $X=(X_t)_{t\geqslant 0}$ and $Y=(Y_t)_{t\geqslant 0}$ be two continuous-time c\'{a}dl\'{a}g martingales with respect to the filtration $(\F_t)_{t\geqslant 0}$ that satisfies the ``usual hypotheses'' \cite{Pro}. In this case the \emph{martingale paraproduct} is also a stochastic process $\big((X\cdot Y)_t\big)_{t\geqslant 0}$, but now defined via the stochastic integral
\begin{equation}
(X\cdot Y)_t := \int_{0+}^t X_{s-} dY_s.
\label{Con:paraprod}
\end{equation}
Since we are allowed to choose dense subspaces on which the initial definition makes sense (and later extend by continuity), we can conveniently assume that $X$ is bounded in $\textup{L}^\infty$ and $Y$ is bounded in $\textup{L}^2$.
We want to prove that \eqref{Con:paraprod} satisfies the same $\textup{L}^p$ estimates as \eqref{Dis:paraprod}.
To do so, we take $(\pi_m)_{m=1}^{\infty}$ to be a refining sequence of partitions
\[ 0 = t_0^{(m)} < t_1^{(m)} < t_2^{(m)} < \dots < t_{n(m)}^{(m)} = t \]
such that $\lim_{m\to\infty}\text{mesh}(\pi_m)=0$. We can calculate \eqref{Con:paraprod} as the limit of the Riemann sums in the following way:
\begin{equation}
\int_{0+}^t X_{s-} dY_s = \lim_{m\to\infty} \sum_{k=1}^{n(m)} X_{t_{k-1}^{(m)}}(Y_{t_k^{(m)}}-Y_{t_{k-1}^{(m)}}).
\label{con:riemann}
\end{equation}
The above limit is interpreted as the convergence in probability; for more details see \cite{Pro}. Notice that the right hand side of \eqref{con:riemann} is actually a limit of discrete-time martingale paraproducts \eqref{Dis:paraprod}. By passing to an a.s.\@ convergent subsequence, using Fatou's lemma, and applying \eqref{dis:est}, we get the desired estimate for \eqref{Con:paraprod}:
\[ \|(X \cdot Y)_t\|_{\textup{L}^{r'}} \leq \sup_{m} \Big\|\sum_{k=1}^{n(m)} X_{t_{k-1}^{(m)}}(Y_{t_k^{(m)}}-Y_{t_{k-1}^{(m)}})\Big\|_{\textup{L}^{r'}}
\lesssim_{p,q,r} \|X_t\|_{\textup{L}^p}\|Y_t\|_{\textup{L}^q} \]
for the exponents $p,q,r$ satisfying \eqref{eq:exponents}.

As a special case we can consider martingales with respect to the augmented filtration of the one-dimensional Brownian motion $(B_t)_{t\geqslant 0}$. For simplicity we also assume that $Y_0=0$, since otherwise we can pass to the martingale $Y_t-Y_0$. Then
\begin{equation}
(X\cdot Y)_t = \int_{0}^{t} X_s dY_s,
\label{BBparaprod}
\end{equation}
because $(X_t)_{t\geqslant 0}$ and $(Y_t)_{t\geqslant 0}$ now a.s.\@ have continuous paths. We remark that \eqref{BBparaprod} are the martingale paraproducts studied by Ba\~{n}uelos and Bennett in \cite{BB:par} and they established $\textup{L}^p$, $\textup{H}^p$, and $\textup{BMO}$ estimates for \eqref{BBparaprod}. Their proof of the $\textup{L}^p$ estimates uses Doob's inequality and the Burkholder-Gundy inequality.

In this particular case we can give yet another short proof, by applying It\={o}'s formula instead of approximating by discrete-time processes. It is more convenient to bound the trilinear form that we obtain by dualizing:
\begin{equation*}
\Lambda_t(X,Y,Z) := \mathbb{E} \big((X\cdot Y)_t Z\big) = \mathbb{E} \big((X\cdot Y)_t (Z_t-Z_0)\big).
\end{equation*}
Here $Z$ is a square-integrable random variable and $Z_s:=\E(Z|\F_s)$ is the corresponding martingale. Using It\={o}'s isometry we get
\[\Lambda_t(X,Y,Z) = \mathbb{E}\bigg( \Big(\int_{0}^{t} X_s dY_s\Big) \Big(\int_{0}^{t} 1 dZ_s\Big) \bigg) = \mathbb{E} \int_{0}^{t} X_s d\langle Y,Z\rangle_s, \]
where $\langle Y,Z\rangle_t$ is the predictable quadratic covariation process, which in the case of the Brownian filtration coincides with the quadratic covariation $[Y,Z]_t$.

Again we assume that $\|X_t\|_{\textup{L}^p}<\infty$, $\|Y_t\|_{\textup{L}^q}<\infty$, $\|Z_t\|_{\textup{L}^r}<\infty$, $X_t,Y_t,Z_t\geqslant 0$, and we introduce three new martingales (for a fixed $t$ and for $s\in[0,t]$):
\[ U_s := \mathbb{E}(X_t^p|\mathcal{F}_s),\quad V_s := \mathbb{E}(Y_t^q|\mathcal{F}_s),\quad W_s := \mathbb{E}(Z_t^r|\mathcal{F}_s).  \]
If we denote $\mathbf{X}_s = (X_s,Y_s,Z_s,U_s,V_s,W_s) = (X_s^i)_{i=1}^6$, then It\={o}'s formula gives us
\[ \mathcal{B}(\mathbf{X}_t) - \mathcal{B}(\mathbf{X}_0) = \sum_{i=1}^6 \int_{0}^{t} \partial_i \mathcal{B}(\mathbf{X}_s) dX_s^i+ \frac{1}{2} \sum_{i,j=1}^6 \int_{0}^{t} \partial_i \partial_j \mathcal{B}(\mathbf{X}_s) d\langle X^i,X^j\rangle_s, \]
where $\B$ is the previously constructed Bellman function. The first term on the right hand side is the martingale part, so by taking the expectation of the above expression, we get
\[\mathbb{E} \big( \mathcal{B}(\mathbf{X}_t) - \mathcal{B}(\mathbf{X}_0) \big) = \mathbb{E} \Big( \int_{0}^{t} \frac{1}{2} \sum_{i,j=1}^6 \partial_i \partial_j \mathcal{B}(\mathbf{X}_s) d\langle X^i,X^j\rangle_s \Big). \]
Using the martingale representation theorem we can write $\mathbf{X}_t$ in the form
\[ \mathbf{X}_t = \mathbf{X}_0 + \int_{0}^{t} \mathbf{A}_s dB_t, \]
where $(\mathbf{A}_t)_{t\geqslant 0}$ is a predictable process. Since $d\langle X^i,X^j\rangle_s = A_s^i A_s^j ds$, using \eqref{eq:b2} and \eqref{eq:b3inf} gives us
\begin{equation*}
\pm \mathbb{E} \int_{0}^{t} X_s^1 A_s^2 A_s^3 ds
\leqslant \mathcal{C}_{p,q,r}\Big(\frac 1p  \|X_t\|_{\textup{L}^p}^p+\frac 1q \|Y_t\|_{\textup{L}^q}^q+\frac 1r \|Z_t\|_{\textup{L}^r}^r \Big ).
\end{equation*}
Finally, $X_s^1 A_s^2 A_s^3 ds=X_s d\langle Y,Z\rangle_s$ and homogenization of the above expression give us the desired estimate.

However, we should emphasize that in order to be able to use It\={o}'s formula, our Bellman function should be of class $\textup{C}^2$ on the whole domain. This is achieved by shrinking the domain slightly and passing to $\mathcal{B}_{\varepsilon}$ as in the next section; we omit the details.

\subsection{Heat flow paraproducts}
In order to be able to use the constructed Bellman function in relationship with the heat equation we should first ``smoothen it up''. Let us fix a non-negative even $\textup{C}^\infty$ function $\varphi$ supported in $(-1,1)^3$ with integral $1$. For any $\varepsilon>0$ we define the function $\mathcal{A}_{\varepsilon}\colon(\varepsilon,\infty)^3\to\mathbb{R}$ by the formula
\[
\mathcal{A}_{\varepsilon}(u,v,w) := \int_{(-\varepsilon,\varepsilon)^3} \varepsilon^{-3} \varphi(\varepsilon^{-1}a, \varepsilon^{-1}b, \varepsilon^{-1}c) \mathcal{A}(u-a,v-b,w-c) da db dc.
\]
In words, $\mathcal{A}_{\varepsilon}$ is the convolution of $\mathcal{A}$ with the $\textup{L}^1$-normalized dilate of $\varphi$. The newly obtained function is clearly of class $\textup{C}^\infty$. We integrate \eqref{eq:a3} translated by $(a,b,c)$ and multiplied by $\varepsilon^{-3} \varphi(\varepsilon^{-1}a, \varepsilon^{-1}b, \varepsilon^{-1}c)$, and then ``symmetrize'' in $(a,b,c)$ and use the fact that $\varphi$ is even. That way we conclude that $\mathcal{A}_{\varepsilon}$ still satisfies the condition \eqref{eq:a3} and consequently also \eqref{eq:a3inf} at every point of its domain. By the formula \eqref{eq:bellman} with $\mathcal{A}_{\varepsilon}$ in the place of $\mathcal{A}$ we can define a $\textup{C}^\infty$ function $\mathcal{B}_{\varepsilon}$ satisfying property \eqref{eq:b3inf} for any $u,v,w>\varepsilon$ and $U\geqslant u^p$, $V\geqslant v^q$, $W\geqslant w^r$. Moreover, property \eqref{eq:a2} is retained up to an additional loss by the factor $\max\{2^p,2^q,2^r\}$, which in turn guarantees \eqref{eq:b2} for some (sufficiently large) constant $\mathcal{C}_{p,q,r}$ independent of $\varepsilon$.

Now suppose that $f,g,h$ are compactly supported $\textup{C}^\infty$ functions on $\mathbb{R}$. Also, let $k(x,t):=\frac{1}{\sqrt{2\pi t}}e^{-\frac{x^2}{2t}}$ be the heat kernel on the real line and $u$ be the heat extension of $f$:
\[ u(x,t) := \int_{\mathbb{R}} f(y) k(x-y,t) dy.\]
Note that $u$ is the solution of the heat equation
$\partial_t u=\frac{1}{2}\partial_x^2 u$
with the initial condition $\lim_{t\to 0+}u(x,t)=f(x)$.
Analogously we define $v$ and $w$ to be the heat extensions of $g$ and $h$.

We can define the \emph{heat paraproduct}, i.e.\@ the paraproduct with respect to the heat semigroup as a trilinear form
\begin{equation}
\Lambda(f,g,h) := \int_{\mathbb{R}} \int_{0}^{\infty} u(x,t) \,\partial_x v(x,t) \,\partial_x w(x,t) \,dt \,dx.
\label{heat:paraprod}
\end{equation}
If we denote
\[ \varphi_s(x) := k(x,s^2),\quad \psi_s(x) := -2^{1/2} s \,\partial_x k(x,s^2) \]
and substitute $t=s^2$, we get a more familiar expression:
\begin{equation}
\Lambda(f,g,h) = \int_{\mathbb{R}} \int_{0}^{\infty} (f\ast\varphi_s)(x) \,(g\ast\psi_s)(x) \,(h\ast\psi_s)(x) \,\frac{ds}{s} \,dx.
\label{heat:para2}
\end{equation}
Smooth paraproducts like \eqref{heat:para2} appear naturally in the proof of the T1 theorem (see \cite{DJ}), although one usually needs to be more flexible when choosing a bump function $\varphi_s$ and a mean zero bump function $\psi_s$.

Again, we want to prove some $\textup{L}^p$ estimates for \eqref{heat:paraprod}, i.e.
\begin{equation*}
|\Lambda(f,g,h)| \lesssim_{p,q,r} \|f\|_{\textup{L}^p(\R)}\|g\|_{\textup{L}^q(\R)}\|h\|_{\textup{L}^r(\R)},
\end{equation*}
where $p,q,r$ are exponents satisfying \eqref{eq:exponents}.
To do so we will imitate the ``heating'' technique by Nazarov and Volberg \cite{NV} or Petermichl and Volberg \cite{PetVol}.

Assume that $f,g,h$ are non-negative and that none of them is identically $0$.
Fix $R>0$, $\delta>0$, $T>2\delta$, and observe that $u(x,t),v(x,t),w(x,t) > \varepsilon$ whenever $x\in[-R,R]$, $t\in[\delta,T-\delta]$ for some sufficiently small $\varepsilon>0$ depending on $R,\delta,T$, and the functions $f,g,h$. We introduce $U,V,W$ as the heat extensions of $f^p,g^q,h^r$ respectively and define
\[ b(x,t) := \mathcal{B}_\varepsilon\big(u(x,t),v(x,t),w(x,t),U(x,t),V(x,t),W(x,t)\big), \]
where $\B_\varepsilon$ is as above.
It is easy to calculate that
\begin{align*}
\big(\partial_t - {\textstyle\frac{1}{2}}\partial_x^2\big) b(x,t)
& = (\nabla\mathcal{B}_\varepsilon)(u,v,w, U, V, W) \cdot \big(\partial_t - {\textstyle\frac{1}{2}}\partial_x^2\big)(u,v,w, U, V, W) \\
& - {\textstyle\frac{1}{2}} (d^2\mathcal{B}_\varepsilon)(u,v,w, U, V, W)(\partial_x u,\partial_x v,\partial_x w, \partial_x U, \partial_x V, \partial_x W).
\end{align*}
(We have omitted writing the variables $x,t$ on the right hand side.)
Since $u,v,w,U,V,W$ all satisfy the heat equation, the first term on the right hand side is zero and by \eqref{eq:b3inf} we get
\[ \big(\partial_t - {\textstyle\frac{1}{2}}\partial_x^2\big) b(x,t) \geqslant  \pm u(x,t) \,\partial_x v(x,t) \,\partial_x w(x,t). \]
It remains to integrate this inequality over $[-R,R]\times[\delta,T-\delta]$ with an appropriate weight, use Green's formula, and then let $\delta\to0$, $R,T\to\infty$. We omit the details and refer to \cite{NV},\cite{PetVol}.

Let us emphasize once again that the previous trick of ``smoothing'' the Bellman function was already used in \cite{NV} and \cite{PetVol} and no explicit formula is needed for its application.


\section*{Acknowledgments}
This work has been supported by the Croatian Science Foundation under the project 3526. We would like to thank the anonymous referee for several useful comments and suggestions that improved the readability of this paper.



\begin{thebibliography}{99}
\bibitem{as:hmf}
M. Abramowitz, I. A. Stegun (Eds.),
\textit{Handbook of mathematical functions with formulas, graphs, and mathematical tables},
Dover Publications, Inc., New York, 1992.

\bibitem{BB:par}
R. Ba\~{n}uelos, A.\ G. Bennett,
\textit{Paraproducts and commutators of martingale transforms},
Proc. Amer. Math. Soc. \textbf{103} (1988), no. 4, 1226--1234.

\bibitem{BO}
R. Ba\~{n}uelos, A. Os\c{e}kowski,
\textit{On the Bellman function of Nazarov, Treil and Volberg},
Math. Z. {\bf 278} (2014), no. 1--2, 385--399.

\bibitem{BMN}
\'{A}. B\'{e}nyi, D. Maldonado, V. Naibo,
\textit{What is a Paraproduct?},
Notices of the Amer. Math. Soc. \textbf{57} (2010), no. 7, 858--860.

\bibitem{Burk1}
D. L. Burkholder,
\textit{Martingale transforms},
Ann. Math. Statist. \textbf{37} (1966), 1494--1504.

\bibitem{Burk2}
D. L. Burkholder,
\textit{Boundary value problems and sharp inequalites for martingale transforms},
Ann. Probab. {\bf 14} (1984), no. 3, 647--702.

\bibitem{CD1}
A. Carbonaro, O. Dragi\v{c}evi\'{c},
\textit{Bellman function and linear dimension-free estimates in a theorem of Bakry},
J. Funct. Anal. {\bf 265} (2013), no. 7, 1085--1104.

\bibitem{CD2}
A. Carbonaro, O. Dragi\v{c}evi\'{c},
\textit{Functional calculus for generators of symmetric contraction semigroups} (2013), Duke Math. J. {\bf 166} (2017), no. 5, 937--974.

\bibitem{Cohn}
D. L. Cohn,
\textit{Measure Theory},
2nd ed., Birkh\"{a}user, 2013.

\bibitem{Davis}
B. Davis,
\textit{On the $L^p$ norms of stochastic integrals and other martingales},
Duke Math. J. {\bf 43} (1976), no. 4, 697--704.

\bibitem{DJ}
G. David, J.-L. Journ\'{e},
\textit{A boundedness criterion for generalized Calder\'{o}n-Zygmund operators},
Ann. of Math. (2) {\bf 120} (1984), no. 2, 371--397.

\bibitem{pd:L4}
P. Durcik,
\textit{An $L^4$ estimate for a singular entangled quadrilinear form},
Math. Res. Lett. {\bf 22} (2015), no. 5, 1317--1332.

\bibitem{pd:Lp}
P. Durcik,
\textit{$L^p$ estimates for a singular entangled quadrilinear form} (2015),
to appear in Trans. Amer. Math. Soc., available at arXiv:1506.08150.

\bibitem{DKST}
P. Durcik, V. Kova\v{c}, K. A. \v{S}kreb, C. Thiele,
\textit{Norm-variation of ergodic averages with respect to two commuting transformations} (2016), to appear in Ergodic Theory Dynam. Systems, available at arXiv:1603.00631.

\bibitem{JP:par}
S. Janson, J. Peetre,
\textit{Paracommutators-Boundedness and Schatten-Von Neumann Properties},
Trans. Amer. Math. Soc. \textbf{305} (1988), no. 2, 467--504.

\bibitem{vk:tp}
V. Kova\v{c},
\textit{Boundedness of the twisted paraproduct},
Rev. Mat. Iberoam. {\bf 28} (2012), no. 4, 1143--1164.

\bibitem{ks:mart}
V. Kova\v{c}, K. A. \v{S}kreb,
\textit{One modification of the martingale transform and its applications to paraproducts and stochastic integrals},
J. Math. Anal. Appl. {\bf 426} (2015), no. 2, 1143--1163.

\bibitem{NT}
F. L. Nazarov, S. R. Treil,
\textit{The hunt for a Bellman function: applications to estimates for singular integral operators
and to other classical problems of harmonic analysis} (in Russian),
Algebra i Analiz \textbf{8} (1996), no. 5, 32--162,
English transl. in St. Petersburg Math. J. \textbf{8} (1997), no. 5, 721--824.

\bibitem{NV}
F. L. Nazarov, A. Volberg,
\textit{Heating of the Ahlfors-Beurling operator and estimates of its norm},
St. Petersburg Math. J. \textbf{15} (2004), no. 4, 563--573.

\bibitem{Osek}
A. Os\c{e}kowski,
\textit{Sharp martingale and semimartingale inequalities},
Monografie Matematyczne 72. Springer, Basel, 2012.

\bibitem{PetVol}
S. Petermichl, A. Volberg,
\textit{Heating of the Ahlfors--Beurling operator: weakly quasiregular maps on the plane are quasiregular},
Duke Math. J. \textbf{112} (2002), no. 2, 281--305.

\bibitem{Pro}
P. E. Protter,
\textit{Stochastic Integration and Differential Equations},
2nd ed., ver. 2.1,
Stoch. Model. Appl. Probab., vol. 21, Springer-Verlag, Berlin, 2005.

\bibitem{Thie}
C. Thiele,
\textit{Wave Packet Analysis},
CBMS Reg. Conf. Ser. Math., \textbf{105}, AMS, Providence, RI, 2006.

\bibitem{VV}
V. Vasyunin, A. Volberg,
\textit{Bellster and others} (2008),
preprint.

\bibitem{w:m9}
Wolfram Research, Inc., \emph{Mathematica}, Ver. 9.0, Champaign, IL, 2012.

\end{thebibliography}
\end{document}